\documentclass[preprint,12pt]{elsarticle}

\usepackage{amsmath,amssymb,amsthm}
\usepackage{hyperref}
\usepackage[a4paper, margin=1in]{geometry}
\usepackage{graphicx}
\usepackage{lmodern}
\usepackage{subcaption}
\usepackage{tikz}
\usetikzlibrary{calc}
\newcommand{\Aut}{\operatorname{Aut}}
\newcommand{\APP}{\operatorname{APP}}

\newtheorem{theorem}{Theorem}[section]
\newtheorem{proposition}[theorem]{Proposition}

\newtheorem{corollary}[theorem]{Corollary}

\theoremstyle{definition}
\newtheorem{definition}[theorem]{Definition}
\newtheorem{example}[theorem]{Example}

\theoremstyle{remark}
\newtheorem{remark}[theorem]{Remark}

\begin{document}

\begin{frontmatter}

\title{Constructing Cospectral Hypergraphs Using Adjacency Preserving Permutations }

\author[inst1]{Arpita Das}
\ead{arpita.das1201@gmail.com}

\author[inst2]{Debjit Ghoshal}
\ead{debjitghoshal181@gmail.com}

\affiliation[inst1]{
  organization={Bethuadahari College},
  city={Bethuadahari},
  country={India}
}

\affiliation[inst2]{
  organization={TIFR Centre for Applicable Mathematics},
  city={Bengaluru},
  country={India}
}

\begin{abstract}
    In this paper, we use permutations of the vertices of hypergraphs that preserve their adjacency matrices defined in \cite{banerjee2021spectrum}. These permutations form a group which lies between the automorphism group and the full symmetric group $S_n$ where $n$ is the number of vertices. Using elements from that group which are not automorphisms, we construct cospectral non-isomorphic hypergraphs. We use hypergraph structure coming from finite projective planes to show existence of infinitely many such hypergraphs. Also this method is very specific to uniform hypergraphs with uniformity $\ge 3$, i.e not applicable to graphs.
\end{abstract}


\end{frontmatter}

\section{Introduction}
\noindent
A hypergraph $\mathcal{H}(V,E)$ consists of a vertex set $V(\mathcal{H})$ and a hyperedge set $E(\mathcal{H})$. Each hyperedge can be formed with more than two vertices and the number of vertices in one such edge is the cardinality of that hyperedge. A hypergraph is said to be $k$-uniform if all hyperedges have cardinality $k$. From a hypergraph one can construct its adjacency matrix, and the multiset of adjacency eigenvalues, counted with algebraic multiplicity, is called the adjacency spectrum or simply the spectrum. Two hypergraphs are called adjacency-cospectral if they have the same adjacency-spectrum. In recent years, spectral hypergraph theory has attracted a lot of attentions \cite{banerjee2021spectrum,cooper2012spectra,feng1996spectra}. In this area one of the interesting problems is the construction of adjacency-cospectral but non-isomorphic hypergraphs. Although the adjacency tensor gives a standard representation of a uniform hypergraph \cite{cooper2012spectra}, extracting spectral information from higher-order tensors can be computationally challenging, with several general tensor eigenvalue problems known to be NP-hard \cite{hillar2013tensorproblemsnphard}. In this paper, we therefore consider the adjacency matrix defined by Banerjee \cite{banerjee2021spectrum}. In contrast to the extensive literature on cospectral graphs, relatively few general methods are known for constructing cospectral hypergraphs. Existing approaches are mainly based on switching techniques. Bu, Zhou, and Wei \cite{bu2014cospectral} extended the classical Godsil--McKay switching to construct $E$-cospectral uniform hypergraphs, while Sarkar and Banerjee \cite{sarkar2020joins} developed a similar switching method for hypergraphs. More recently, Abiad and Khramova \cite{monge2025constructing} proposed new switching-based methods for constructing uniform hypergraphs that are cospectral with respect to their adjacency matrices and adjacency tensors. In this paper, we introduce an approach for constructing adjacency-cospectral but non-isomorphic uniform hypergraphs. Our construction is based on the notion of adjacency-preserving permutations introduced by Mulas, Sánchez-García, and MacArthur \cite{Mulas_2021}. Though it was named adjacency-automorphism there, we will call the collection of such permutations as $\APP(\mathcal{G})$, i.e. adjacency preserving permutations for a hypergraph $\mathcal{G}$. Starting from a uniform hypergraph that has a nontrivial adjacency-preserving permutation which is not an automorphism, we construct non-isomorphic, adjacency-cospectral hypergraphs by attaching pendant edges—edges in which all but one vertex have degree one. Our method exploits the symmetry in adjacency-matrix level induced by adjacency-preserving permutations.

\section{Adjacency Preserving Permutation Group}
\noindent
We know that for simple graphs, their adjacency matrices uniquely determines the graph. For hypergraph this may not be the case. In fact, two different hypergraphs may have same adjacency matrix. This motivates the notion of coadjacent hypergraphs and adjacency preserving permutations introduced in this section. 
\begin{definition}[Coadjacent Hypergraphs]{\label{coadjacen-defn}}
Let $\mathcal{G}=(V,E)$ and $\mathcal{H}=(V,F)$ be hypergraphs on the same vertex set
$V=\{1,2,\ldots,n\}$. We say that $\mathcal{G}$ and $\mathcal{H}$ are
\emph{coadjacent} if $A_{\mathcal G}=A_{\mathcal H}$ (with a fixed ordering of vertices),
where the adjacency matrix of a hypergraph is defined by \cite{banerjee2021spectrum}

\begin{equation}
    (A_{\mathcal G})_{ij}
=
\begin{cases}
\displaystyle
\sum_{\substack{e \in E \\ i,j \in e}}
\frac{1}{|e|-1},
& \text{if } i \neq j, \\[1.2ex]
0,
& \text{if } i=j.
\end{cases}
\label{adj-defn}
\end{equation}

\end{definition}

\begin{example}\label{ex:coadjacent}
Consider the vertex set
\[
V = \{1,2,3,4,5,6\},
\]
and the \(3\)-uniform hypergraphs
\[
\mathcal{G}_1 = (V,E_1)
\qquad \text{and} \qquad
\mathcal{G}_2 = (V,E_2),
\]
where
\begin{align*}
E_1 &= \bigl\{
\{1,2,3\},\,
\{3,4,5\},\,
\{5,6,1\},\,
\{2,4,6\}
\bigr\},\\[6pt]
E_2 &= \bigl\{
\{6,1,2\},\,
\{2,3,4\},\,
\{4,5,6\},\,
\{1,3,5\}
\bigr\}.
\end{align*}

\begin{figure}[ht]
\centering
\begin{subfigure}[b]{0.25\textwidth}
\centering
\includegraphics[width=\textwidth]{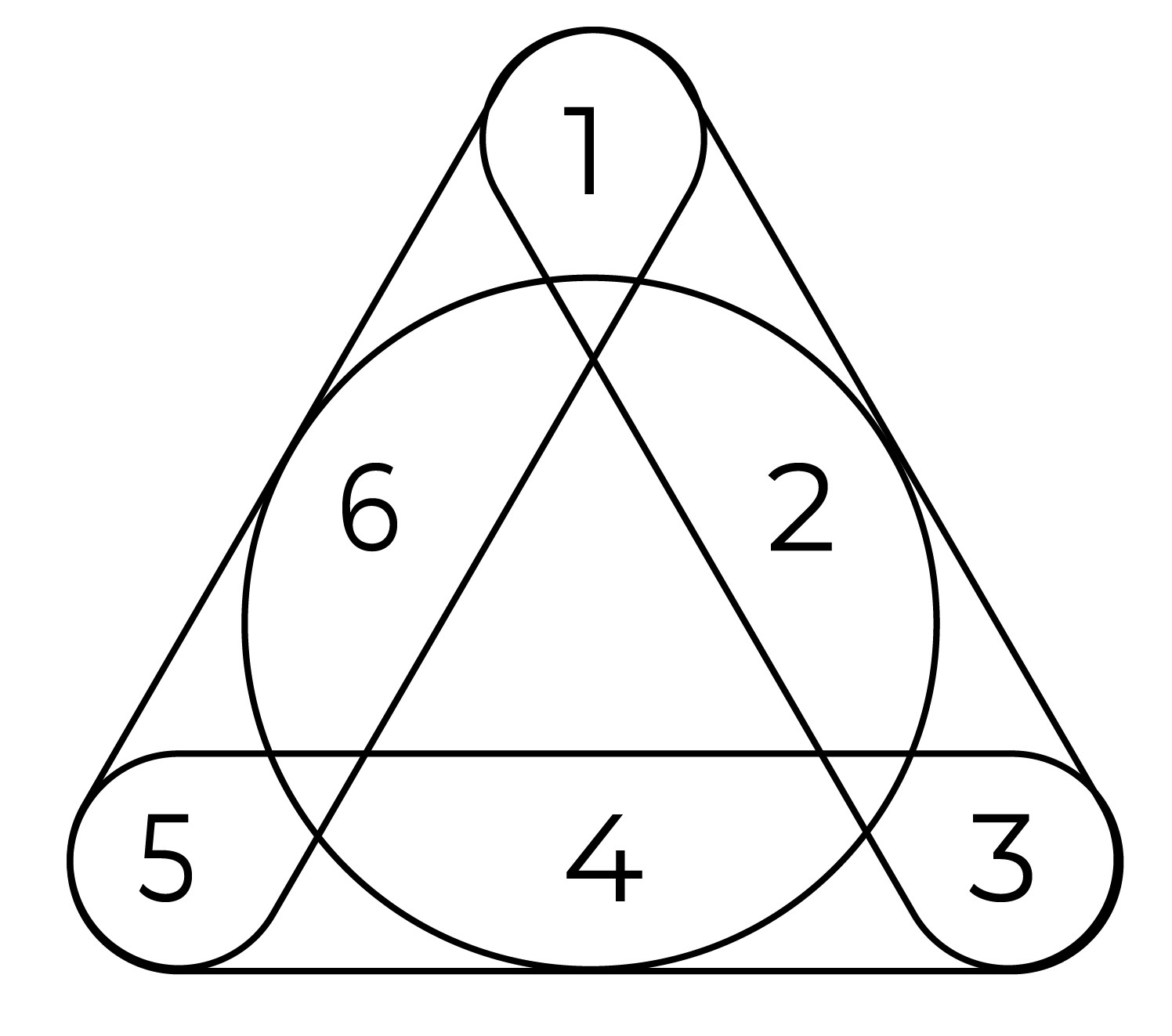}
\caption{$\mathcal G_1$}
\end{subfigure}
\hspace{0.2\textwidth}
\begin{subfigure}[b]{0.25\textwidth}
\centering
\includegraphics[width=\textwidth]{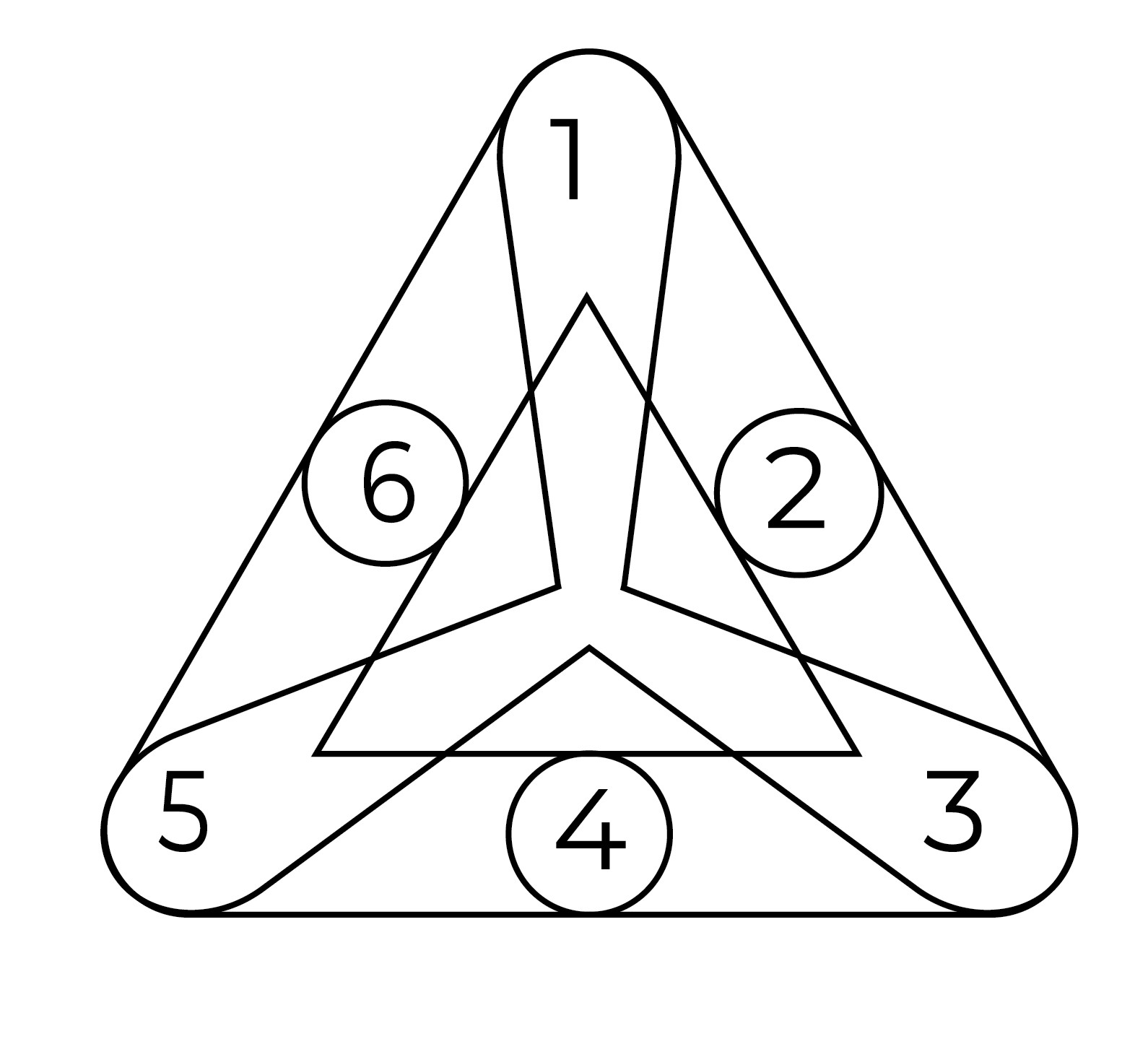}
\caption{$\mathcal G_2$}
\end{subfigure}
\caption{Two coadjacent hypergraphs.}
\label{fig:coadjacent}
\end{figure}
\noindent
Using \eqref{adj-defn}, direct computation gives
\[
A_{\mathcal G_1}
=
A_{\mathcal G_2}
=
\begin{bmatrix}
0 & \tfrac12 & \tfrac12 & 0 & \tfrac12 & \tfrac12\\
\tfrac12&0&\tfrac12&\tfrac12&0&\tfrac12\\
\tfrac12&\tfrac12&0&\tfrac12&\tfrac12&0\\
0&\tfrac12&\tfrac12&0&\tfrac12&\tfrac12\\
\tfrac12&0&\tfrac12&\tfrac12&0&\tfrac12\\
\tfrac12&\tfrac12&0&\tfrac12&\tfrac12&0
\end{bmatrix}
=
\frac12\,J_2\otimes(J_3-I_3).
\]
\noindent
Hence $\mathcal G_1$ and $\mathcal G_2$ are coadjacent. Moreover, they are even
isomorphic, since the permutation $\pi=(1\,2\,3\,4\,5\,6)$
maps $E_1$ onto $E_2$.
\end{example}
\noindent
The previous example shows that isomorphic hypergraphs may be coadjacent. The next example shows that the converse can be false.

\begin{example}\label{ex:coadjacent-nonisomorphic}
Let
\[
V = \{1,2,3,4,5,6,7\}
\]
and consider the hypergraphs
\[
\mathcal{H}_1 = (V,E_1)
\qquad \text{and} \qquad
\mathcal{H}_2 = (V,E_2),
\]
where
\begin{align*}
E_1 = \Bigl\{
&\{1,2,3\},\,
\{3,4,5\},\,
\{5,6,1\},\,
\{2,4,6\},\,
\{3,4,7\},\,
\{4,5,7\}
\Bigr\},
\\[6pt]
E_2 = \Bigl\{
&\{6,1,2\},\,
\{2,3,4\},\,
\{4,5,6\},\,
\{1,3,5\},\,
\{3,4,7\},\,
\{4,5,7\}
\Bigr\}.
\end{align*}
\begin{figure}[ht]
    \centering
    \begin{subfigure}[b]{0.25\textwidth}
        \centering
        \includegraphics[width=\textwidth]{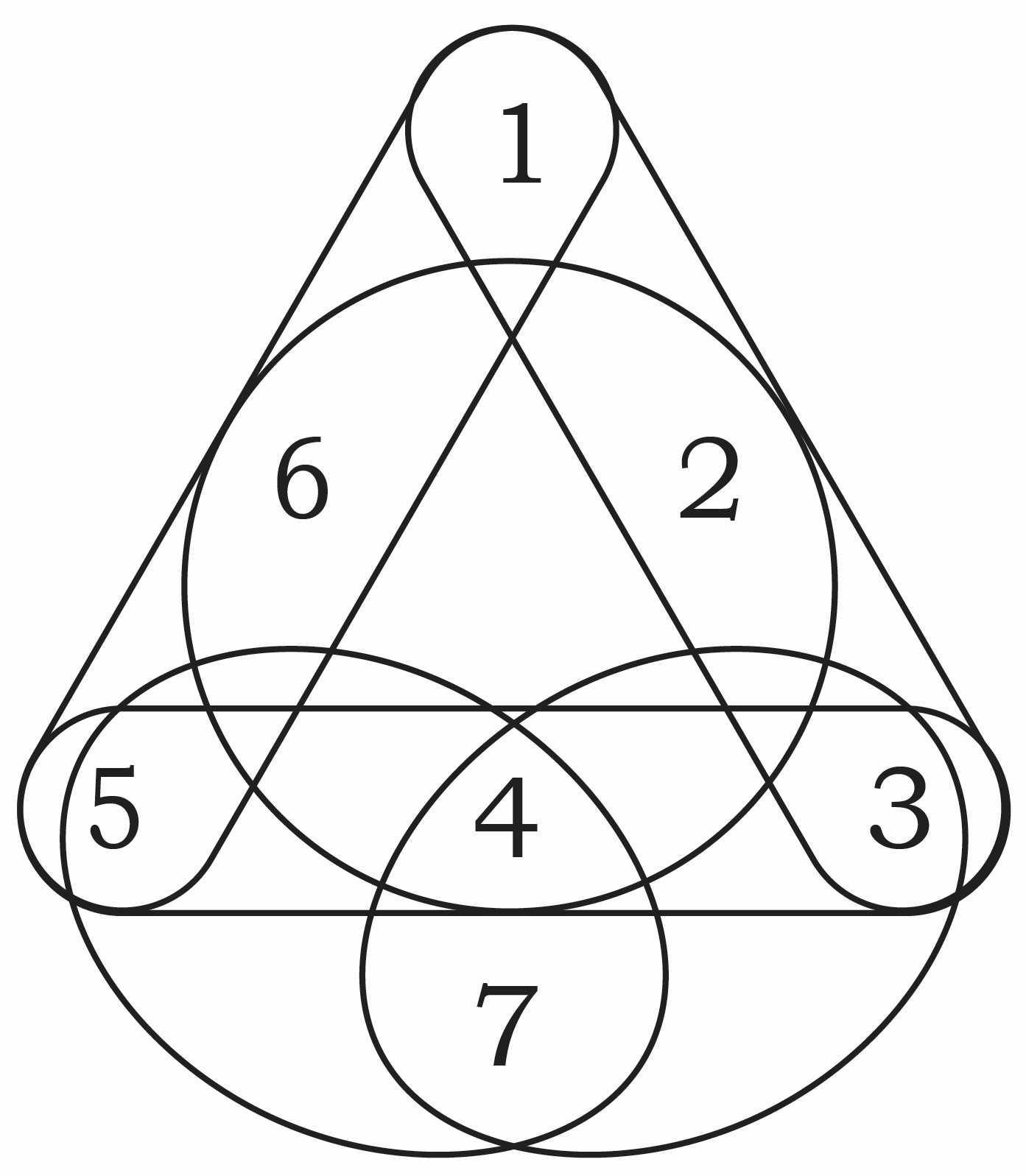}
        \caption{$\mathcal H_1$}
    \end{subfigure}
   \hspace{0.2\textwidth}
    \begin{subfigure}[b]{0.25\textwidth}
        \centering
        \includegraphics[width=\textwidth]{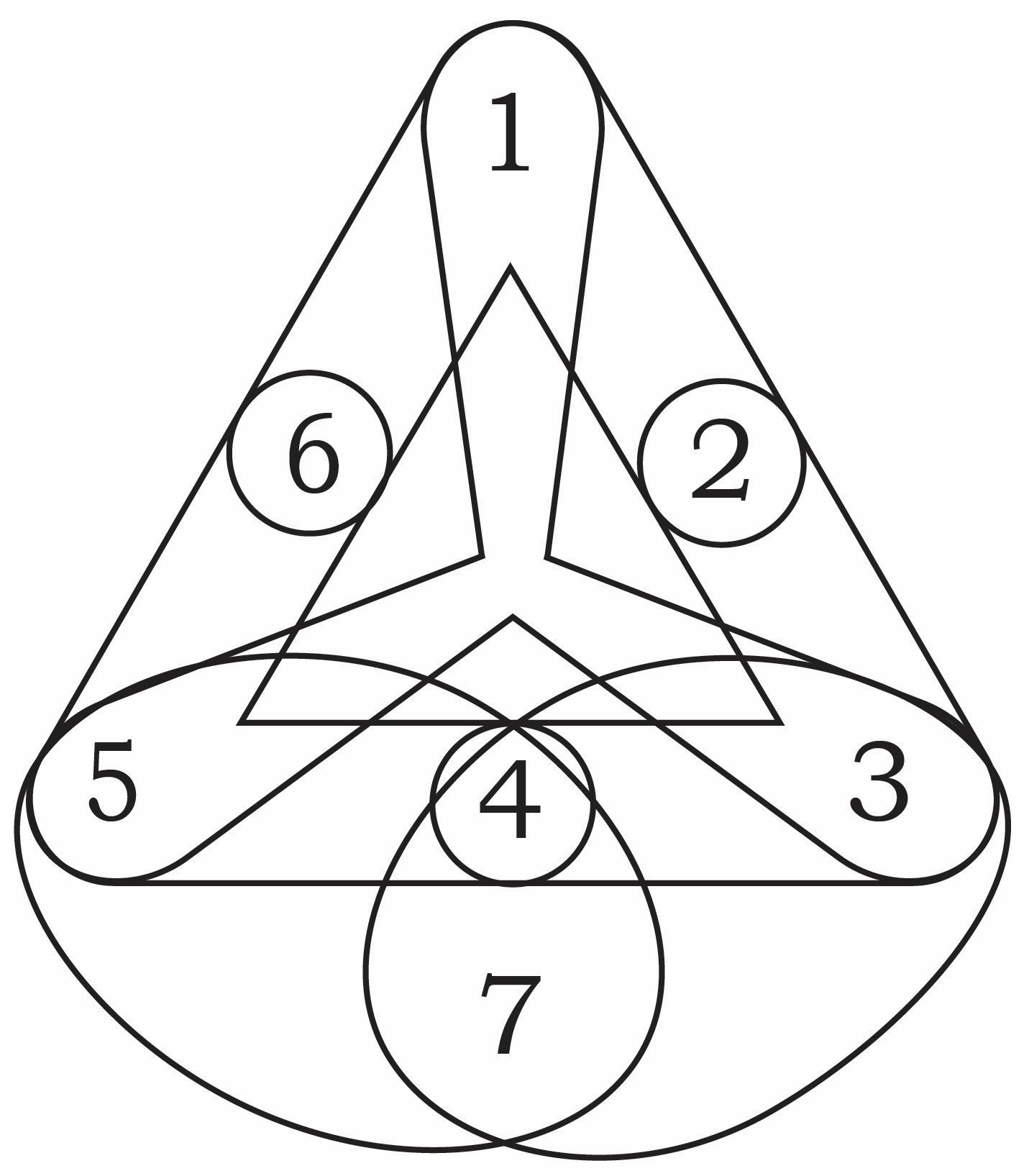}
        \caption{$\mathcal H_2$}
    \end{subfigure}
    \caption{Two coadjacent non-isomorphic hypergraphs.}
    \label{fig:nonisoorphic}
\end{figure}

\noindent
Clearly $\mathcal H_1$ and $\mathcal H_2$ are coadjacent as their adjacency matrices coincide. Indeed,
\[
A_{\mathcal H_1}
=
A_{\mathcal H_2}
=
\begin{bmatrix}
0 & \tfrac12 & \tfrac12 & 0 & \tfrac12 & \tfrac12 & 0 \\
\tfrac12 & 0 & \tfrac12 & \tfrac12 & 0 & \tfrac12 & 0 \\
\tfrac12 & \tfrac12 & 0 & 1 & \tfrac12 & 0 & \tfrac12 \\
0 & \tfrac12 & 1 & 0 & 1 & \tfrac12 & 1 \\
\tfrac12 & 0 & \tfrac12 & 1 & 0 & \tfrac12 & \tfrac12 \\
\tfrac12 & \tfrac12 & 0 & \tfrac12 & \tfrac12 & 0 & 0 \\
0 & 0 & \tfrac12 & 1 & \tfrac12 & 0 & 0
\end{bmatrix}.
\]

\noindent
To see that \(\mathcal{H}_1\) and \(\mathcal{H}_2\) are non-isomorphic, observe in \(\mathcal{H}_1\), the edge $\{3,4,5\}$ contains one degree-\(4\) vertex and two degree-\(3\) vertices. However, \(\mathcal{H}_2\) has no edge with this property. Observe that coadjacency will surely give cospectrality , but the converse is not true in general.
\end{example}
\begin{remark}
For $2$-uniform hypergraphs, that is, ordinary graphs, coadjacency coincides with graph equality. Indeed, if two graphs have the same adjacency matrix with respect to a fixed vertex ordering, then they have exactly the same edge set. Thus, the phenomenon illustrated in Example~\eqref{ex:coadjacent-nonisomorphic} is specific to hypergraphs which has alteast one edge with cardinality $\ge 3$.
\end{remark}

\noindent
Example~\eqref{ex:coadjacent}  and ~\eqref{ex:coadjacent-nonisomorphic} shows that a relabeling of the vertices may leave the adjacency matrix unchanged when expressed with respect to a fixed ordering of the vertex set. This observation motivates the following notion, adapted from the notion of an adjacency automorphism for oriented hypergraphs studied by Mulas, Sánchez-García, and MacArthur \cite{Mulas_2021}, which plays a central role in our construction. We use the term \emph{adjacency-preserving permutation} for what is referred to as an \emph{adjacency automorphism} in \cite{Mulas_2021}.
\begin{definition}[Adjacency-Preserving Permutation {\cite{Mulas_2021}}]
\label{def:adj-pres}
Let $\mathcal G=(V,E)$ be a hypergraph with
$V=\{1,2,\ldots,n\}$, and let
$\pi\in S_n$.
Denote by $\mathcal G^\pi=(V,E^\pi)$, the hypergraph obtained from $\mathcal G$ by relabeling every vertex $v$ by $\pi(v)$, where
$E^\pi=\{\pi(e):e\in E\},
\quad
\pi(e)=\{\pi(v):v\in e\}.$ We say that $\pi$ is an
\emph{adjacency-preserving permutation} if
$A_{\mathcal G}=A_{\mathcal G^\pi}$, where both adjacency matrices are computed with respect to the fixed vertex ordering $(1,2,\ldots,n)$.
\end{definition}

\noindent
The previous example immediately gives a nontrivial adjacency-preserving
permutation.
\begin{example}
Consider the example~\eqref{ex:coadjacent}, let $\pi=(1\,2\,3\,4\,5\,6)$. Since $\pi$ maps the edge set of $\mathcal G_1$ onto that of
$\mathcal G_2$, we have $\mathcal G_1^\pi=\mathcal G_2$.
Moreover, $A_{\mathcal G_1}=A_{\mathcal G_2}$. Hence $A_{\mathcal G_1}=A_{\mathcal G_1^\pi}$, so $\pi$ is an adjacency-preserving permutation of
$\mathcal G_1$. Notice that $\pi$ is not an automorphism of $\mathcal G_1$, since
$E_1 ^\pi= E_2 \neq E_1$.
\end{example} 

\noindent

The collection of all adjacency-preserving permutations forms a group.
\begin{definition}[Adjacency-Preserving Permutation Group]
Let $\mathcal G=(V,E)$ be a hypergraph on
$V=\{1,2,\ldots,n\}$.
The
\emph{adjacency-preserving permutation group}
of $\mathcal G$ is
\[
\APP(\mathcal G)
=
\left\{
\pi\in S_n:
A_{\mathcal G}
=
A_{\mathcal G^\pi}
\right\}.
\]
Equivalently, if \(P_\pi\) denotes the permutation matrix corresponding to \(\pi\), then
\[
A_{\mathcal G^\pi}
=
P_\pi A_{\mathcal G}P_\pi^{T}.
\]
Hence, the group of adjacency-preserving permutations can be written as 
\begin{equation}
    \APP(\mathcal G)
=
\left\{
\pi\in S_n:
P_\pi A_{\mathcal G}P_\pi^{T}=A_{\mathcal G}
\right\}.
\label{app-defn}
\end{equation}

\end{definition}

\begin{proposition}
For every hypergraph $\mathcal G$,
$\APP(\mathcal G)$ is a subgroup of
$S_n$.
\end{proposition}

\begin{proof}
Clearly identity permutation on the vertex set belongs to $\APP(\mathcal G)$. So $\APP(\mathcal G)$ is a non-empty subset of the group $S_n$. 
\noindent
Now let $\pi,\sigma\in\APP(\mathcal G)$. It is enough to show that $$\sigma\circ\pi^{-1}\in\APP(\mathcal G).$$ Now using the fact that
$P_{\sigma\circ\pi^{-1}}=P_\sigma P_\pi^{-1}$ and $P_\pi^{-1} =P_\pi^T$ , we have 
\begin{align*}
A_{\mathcal{G}^{\sigma\circ\pi^{-1}}}
&= P_\sigma P_\pi^{-1} A_{\mathcal{G}} (P_\pi^{-1})^T P_\sigma^T \\
&= P_\sigma \left(P_\pi^T A_{\mathcal{G}} P_\pi\right) P_\sigma^T \\
&= P_\sigma A_{\mathcal{G}} P_\sigma^T
\qquad \text{(since } \pi \in \APP(\mathcal{G})\text{)} \\
&= A_{\mathcal{G}} \qquad \text{(since } \sigma \in \APP(\mathcal{G})\text{)}.
\end{align*}

\noindent
So $\sigma\circ\pi^{-1} \in \APP(\mathcal G)$.
Hence it is a subgroup of $S_n$.
\end{proof}

\noindent
The adjacency-preserving permutation group naturally contains the automorphism group.

\begin{proposition}{\label{pr:subgroup}}
Let $\mathcal G=(V,E)$ be a hypergraph. Then $\Aut(\mathcal G) \le \APP(\mathcal G).$
\end{proposition}

\begin{proof}
Let
$\pi\in\Aut(\mathcal G)$. Then $E^\pi=E \implies \mathcal G^\pi=\mathcal G \implies A_{\mathcal G^\pi} = A_{\mathcal G},$ which gives $\pi\in \APP(\mathcal G).$ Hence
$\Aut(\mathcal G) \subseteq \APP(\mathcal G).$ Since both are subgroups of $S_n$, we have $\Aut(\mathcal G) \le \APP(\mathcal G).$
\end{proof}
\noindent
The inclusion above may be strict. Indeed, example~\eqref{ex:coadjacent} gives a permutation that preserves the adjacency matrix but does not preserve the edge set. Hypergraphs admitting such nontrivial adjacency-preserving permutations form the foundation of our construction.

\section{Coadjacent Hypergraphs with Altered Edge Sets}
\noindent
While the automorphism group $\mathrm{Aut}(\mathcal{G})$ consists of all vertex permutations that preserve the edge set exactly, the \emph{adjacency-preserving permutation group} $\mathrm{APP}(\mathcal{G})$ is generally larger, as it includes all permutations that preserve the adjacency matrix—even when the edge structure itself changes. This distinction captures a more relaxed form of symmetry, one that is not constrained by strict edge preservation, but rather by the pairwise co-membership encoded in the adjacency matrix. The set difference $\mathrm{APP}(\mathcal{G}) \setminus \mathrm{Aut}(\mathcal{G})$ thus represents permutations that preserve the adjacency structure altering the edge set.

\medskip

\begin{proposition}\label{aut=app}
Let \( \mathcal{G} = (V, E) \) be a simple graph (i.e., a 2-uniform hypergraph). Then the adjacency-preserving permutation group coincides with the automorphism group, $\mathrm{APP}(\mathcal{G}) = \mathrm{Aut}(\mathcal{G})$.

\end{proposition}

\begin{proof}
Let \(V=\{1,2,\dots,n\}\), and let \(A=A_{\mathcal G}\) be the adjacency matrix of \(\mathcal G\), defined by
\[
A_{ij}=
\begin{cases}
1 & \text{if } \{i,j\}\in E,\\
0 & \text{otherwise}.
\end{cases}
\]
Observe that this definition coincides with the earlier definition in \eqref{coadjacen-defn}.

\noindent
It is enough to show that
\[
\mathrm{APP}(\mathcal G)\subseteq \mathrm{Aut}(\mathcal G),
\]

\noindent
Let \(\pi\in\mathrm{APP}(\mathcal G)\), and consider the graph
\[
\mathcal G^\pi=(V,E^\pi),
\qquad
E^\pi=\{\pi(e):e\in E\}.
\]
By the definition of \(\mathrm{APP}(\mathcal G)\),
\[
A_{\mathcal G^\pi}=A_{\mathcal G}.
\]
For \(i,j\in V\), we have
\[
\begin{aligned}
A_{\mathcal G^\pi}(i,j)=1
&\iff \{i,j\}\in E^\pi\\
&\iff \{\pi^{-1}(i),\pi^{-1}(j)\}\in E.
\end{aligned}
\]
Hence
\[
A_{\mathcal G^\pi}(i,j)
=
A_{\mathcal G}(\pi^{-1}(i),\pi^{-1}(j)).
\]
Since \(A_{\mathcal G^\pi}=A_{\mathcal G}\), it follows that
\[
A_{\mathcal G}(i,j)
=
A_{\mathcal G}(\pi^{-1}(i),\pi^{-1}(j))
\qquad\text{for all }i,j\in V.
\]
Equivalently,
\[
\{i,j\}\in E
\iff
\{\pi^{-1}(i),\pi^{-1}(j)\}\in E.
\]
Replacing \(i,j\) by \(\pi(i),\pi(j)\), we get
\[
\{\pi(i),\pi(j)\}\in E
\iff
\{i,j\}\in E.
\]
Thus \(\pi\) preserves the edge set, and hence
\[
\pi\in\mathrm{Aut}(\mathcal G).
\]
Therefore,
\[
\mathrm{APP}(\mathcal G)\subseteq\mathrm{Aut}(\mathcal G).
\]
Together with the proposition \eqref{pr:subgroup}, we conclude that
\[
\mathrm{APP}(\mathcal G)=\mathrm{Aut}(\mathcal G).
\]
\end{proof}
\noindent
While for simple graphs $\mathcal{G}$ we showed that $\mathrm{Aut}(\mathcal{G}) = \mathrm{APP}(\mathcal{G})$, this is not the case for general hypergraphs.
\begin{example}
\label{strict-inclusion}
    From example \eqref{ex:coadjacent} consider the 3-uniform hypergraph $\mathcal{G}_1 = (V, E_1)$ with
\[
E_1 = \left\{ \{1, 2, 3\}, \{3, 4, 5\}, \{5, 6, 1\}, \{2, 4, 6\} \right\},
\]
The adjacency matrix $A_{\mathcal{G}_1}$ for this hypergraph is given by
\[
A_{\mathcal{G}_1} = \frac12\,J_2\otimes(J_3-I_3).
\]
\noindent
Consider the permutation $\pi = (1\ 2\ 3\ 4\ 5\ 6)$, which cyclically permutes the vertices. Under this permutation, the edge set becomes
\[
E_1^\pi = \left\{ \{6, 1, 2\}, \{2, 3, 4\}, \{4, 5, 6\}, \{1, 3, 5\} \right\},
\]
Which is exactly the hypergraph $\mathcal{G}_2=(V,E_2)$. We know that $A_{\mathcal{G}_1} = A_{\mathcal{G}_2}$. Hence, $\pi \in \mathrm{APP}(\mathcal{G}_1)$. However, $\pi$ does not preserve the edge set, as $E_1 \neq E_2$. Thus, $\pi \notin \mathrm{Aut}(\mathcal{G}_1)$. Therefore, we have $\mathrm{Aut}(\mathcal{G}_1) \subsetneq \mathrm{APP}(\mathcal{G}_1)$.
\end{example}
\noindent
 So adjacency-preserving permutation group can be strictly larger than the automorphism group for higher-order hypergraphs.
To construct cospectral non-isomorphic hypergraphs, we use hypergraphs \( \mathcal{G} \) for which \( \mathrm{Aut}(\mathcal{G}) \) is a proper subset of \( \mathrm{APP}(\mathcal{G}) \). Specifically, we are interested in the non-empty difference of these two groups, i.e., $\mathrm{APP}(\mathcal{G}) \setminus \mathrm{Aut}(\mathcal{G}) \neq \emptyset$. Such hypergraphs serve as the foundation for the construction of cospectral non-isomorphic hypergraphs in this paper.

\section{Constructing Cospectral non-isomorphic hypergraphs}
\noindent
We now arrive at our main result. Our objective is to construct two hypergraphs that are cospectral but non-isomorphic.We use one perturbed edge to introduce new vertices and construct additional edges in such a way that the underlying matrix structure remains unchanged. The construction is as follows.

\noindent
Let $\mathcal{G}=(V,E)$ be a hypergraph with the following three properties
\begin{itemize}
    \item[1.] It is $m$-uniform ($m \ge 3$).
    \item[2.] Every vertex has \textbf{degree atleast $2$}
    \item[3.] $\APP(\mathcal{G})\setminus
\Aut(\mathcal{G})\neq\emptyset.$
\end{itemize}

\noindent
Choose $$\pi\in
\APP(\mathcal{G})\setminus
\Aut(\mathcal{G}).$$

\noindent
Since $\pi$ is not an automorphism of $\mathcal{G}$, $E^\pi\neq E.$ Hence there exists a hyperedge
\[
e=\{v_1,v_2,\ldots,v_m\}\in E
\]
such that
\begin{equation}
    e\notin E^\pi.
    \label{condition 1}
\end{equation}

\noindent
Introduce $m(m-1)$ new vertices
\[
\{a_{i,j}:1\le i\le m,\;1\le j\le m-1\},
\]
all distinct and disjoint from $V$. Define
\[
\widetilde V
=
V\cup
\{a_{i,j}:1\le i\le m,\;1\le j\le m-1\}.
\]

\noindent
For each $i=1,\ldots,m$, define the edge
\[
\widetilde e_i
=
\{v_i,a_{i,1},a_{i,2},\ldots,a_{i,(m-1)}\}.
\]
Now consider the hypergraph
\[
\widetilde{\mathcal G}
=
(\widetilde V,\widetilde E),
\]
with
\[
\widetilde E
=
E\cup
\{\widetilde e_1,\ldots,\widetilde e_m\}.
\]
These extra edges $\widetilde{e_i}$ are called pendant edges as all but one vertices has degree $1$. 
Similarly, define another hypergraph
\[
\widetilde{\mathcal H}
=
(\widetilde V,F),
\]
where
\[
F
=
E^\pi
\cup
\{\widetilde e_1,\ldots,\widetilde e_m\}.
\]
Behind this construction we draw the idea from example \eqref{ex:coadjacent-nonisomorphic}, where we add new vertices and hyperedges in a distinct manner and permute vertices of the original hypergraph. The pendant edges serve as a fixed follower of the corresponding vertex from the original hypergraph.

\begin{theorem}\label{thm:final}
Under the above conditions and constructions we have
\[
\widetilde{\mathcal G}\not\cong
\widetilde{\mathcal H}, \qquad A_{\widetilde{\mathcal G}} = A_{\widetilde{\mathcal H}}
\] 
\end{theorem}

\noindent
Before proceeding to the proof, let us see an example.
\begin{example}\label{finalex}
    Consider the  hypergraph $\mathcal G=(V,E)$ in example \eqref{ex:coadjacent}.
\[
V=\{1,2,3,4,5,6\}, \quad E=
\bigl\{
\{1,2,3\},
\{3,4,5\},
\{5,6,1\},
\{2,4,6\}
\bigr\},
\]
Observe that it is $3$-uniform, each vertex has degree $2$.
Now following example \eqref{strict-inclusion} where we get that the permutation $\pi = (6\ 1\ 2\ 3\ 4\ 5) \in \APP(\mathcal{G})\setminus\Aut(\mathcal{G})$. So it satisfies all the hypotheses we need in Theorem \eqref{thm:final}. 

\medskip
\noindent
Consider the edge $e=\{3,4,5\}$. Clearly, $e \notin E^{\pi}$. 

\medskip
\noindent
Now introduce new vertices $\{7,8,9,10,11,12\}$. 
Consider the hypergraph
$\mathcal{H}_1=(\widetilde{V},\widetilde{E})$ where
\[
\widetilde{V}=\{1,2,3,4,5,6,7,8,9,10,11,12\} = V \cup \{7,8,9,10,11,12\},
\]
\[
\widetilde{E}=
\bigl\{
\{1,2,3\},
\{3,4,5\},
\{5,6,1\},
\{2,4,6\},
\{3,7,8\},
\{4,9,10\},
\{5,11,12\}
\bigr\},
\]
Also consider the hypergraph $\mathcal{H}_2=(\widetilde{V},F)$
where 
\[
\begin{aligned}
F
&= E^{\pi}
   \cup
   \bigl\{
      \{3,7,8\},
      \{4,9,10\},
      \{5,11,12\}
   \bigr\} \\[2mm]
&=
\bigl\{
   \{6,1,2\},
   \{2,3,4\},
   \{4,5,6\},
   \{1,3,5\},
   \{3,7,8\},
   \{4,9,10\},
   \{5,11,12\}
\bigr\}.
\end{aligned}
\]
\begin{figure}[ht]
\centering
\begin{subfigure}[b]{0.25\textwidth}
\centering
\includegraphics[width=\textwidth]{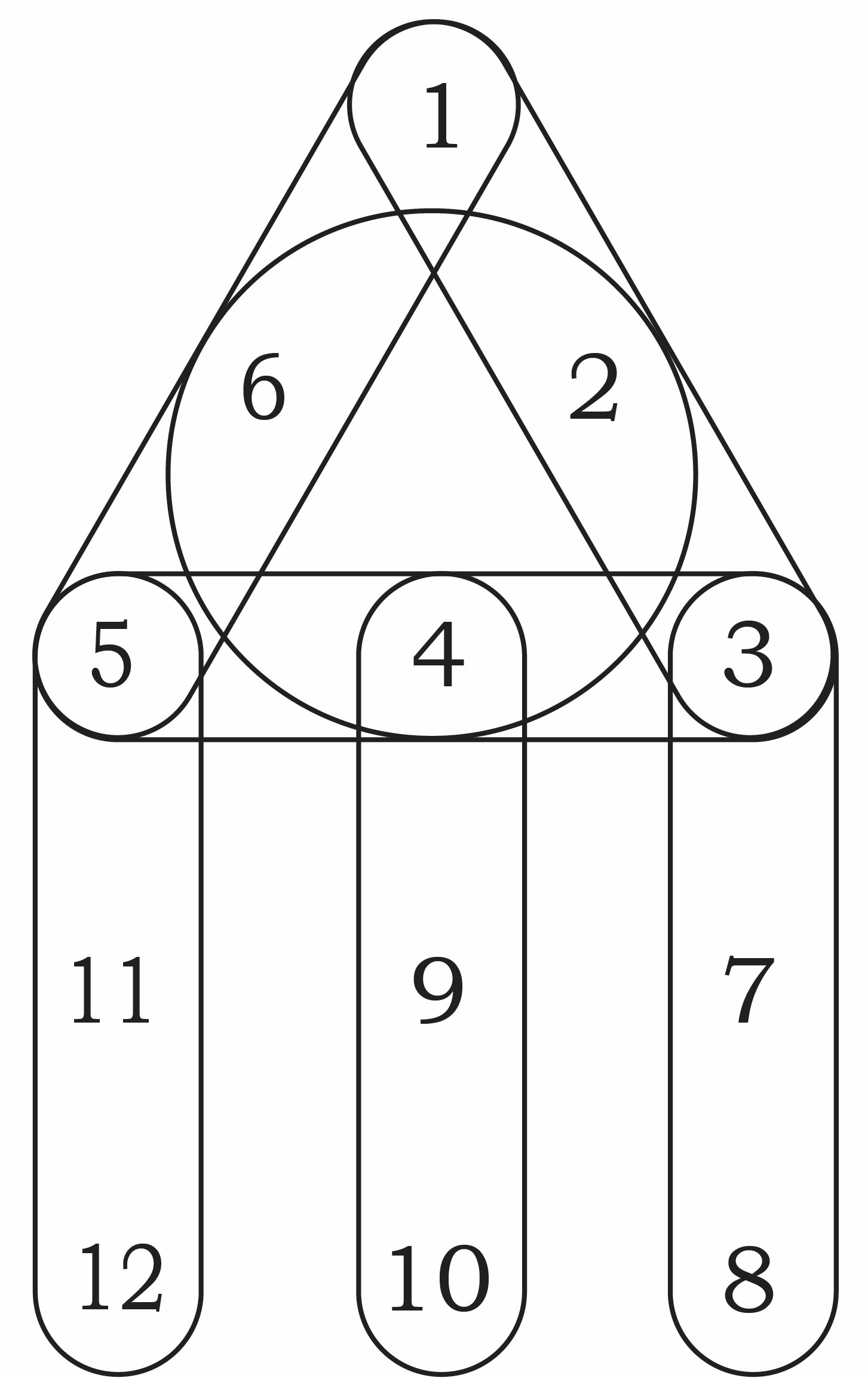}
\caption{$\mathcal H_1$}
\end{subfigure}
\hspace{0.2\textwidth}
\begin{subfigure}[b]{0.25\textwidth}
\centering
\includegraphics[width=\textwidth]{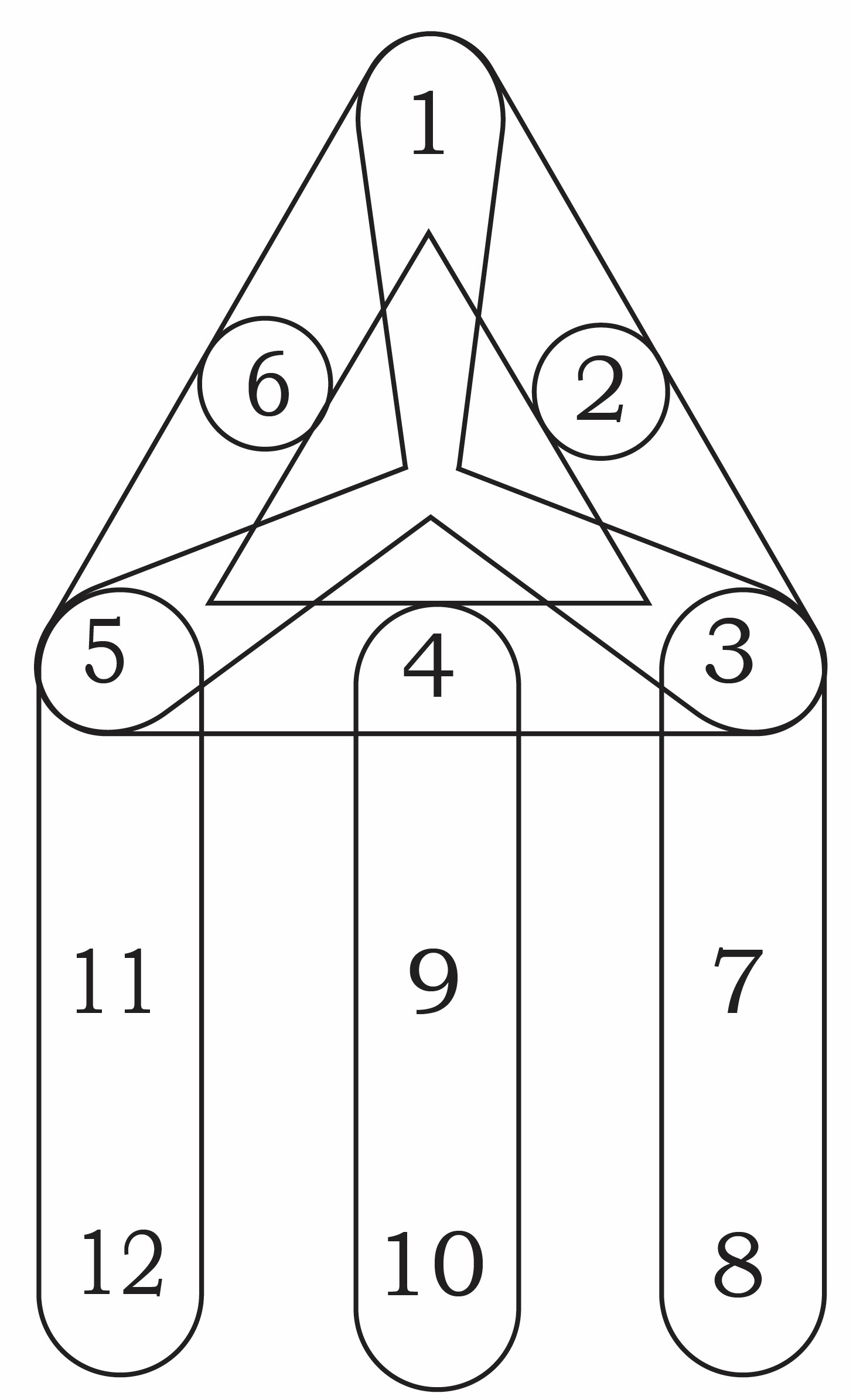}
\caption{$\mathcal H_2$}
\end{subfigure}
\caption{Two coadjacent non-isomorphic hypergraphs.}
\label{fig:nonisomorphic}
\end{figure}
Now proving Theorem \eqref{thm:final} will imply that $\mathcal{H}_1$ and $\mathcal{H}_2$ are non-isomorphic but cospectral.
\end{example}

\begin{proof}
    Let if possible $\Phi$ be an isomorphism between $\widetilde{\mathcal G}$ and $\widetilde{\mathcal H}$. So $\Phi :\widetilde{V}
\longrightarrow
\widetilde{V}$ is a bijection satisfying 
\begin{itemize}
    \item $\Phi$ preserves degree of vertices.
    \item If $u\sim v $ then $\Phi(u) \sim \Phi(v)$.
\end{itemize}
Now let $S=\{a_{i,j}:1\le i\le m,\;1\le j\le m-1\}$ be the set of newly introduced vertices and $D = \{e_1,\cdots,e_m\}$ be the new pendant edges. So $\widetilde{V}=V \cup S$. Every vertex in $S$ has degree one, whereas every vertex in $V$ has degree at least two in both constructed hypergraphs. Since an isomorphism preserves vertex degrees, we have $\Phi(S)=S$. Each vertex in $S$ belongs to exactly one hyperedge, that is its corresponding pendant edge. Since $\Phi$ preserves hyperedges, it maps pendant edges to pendant edges. That means 
$$ \forall e_i  \in D \quad \exists e_j \in D \text{ such that } \Phi(e_i)=e_j, \quad 1\le i,j \le m.$$
As it is an isomorphism, no two $e_i$ can have same image. So we have a permutation $\sigma \in S_m$ such that 
\begin{equation}
\label{key}
    \Phi(e_i)=e_{\sigma(i)}.
\end{equation}
Now observe one crucial thing, each of these $e_i$ has exactly one vertex that is in $V$, i.e  the vertex $v_i$ and it has degree $\ge 2$. From \eqref{key} and using the fact that each pendant edge contains exactly one vertex from $V$, we can say that 
\begin{equation}\label{key1}
    \Phi(v_i) = v_{\sigma(i)}.
\end{equation}
Now observe that $\{v_1,\cdots,v_m\}$ forms an edge $e$ in $\widetilde{\mathcal{G}}$. So $\{\Phi(v_1),\cdots,\Phi(v_m)\}$ also forms an edge in $\widetilde{\mathcal H}.$ Using \eqref{key1}, $\{\Phi(v_1),\cdots,\Phi(v_m)\} = \{v_1,\cdots,v_m\}=e$. That means $\Phi(e)=e \in F$. 

\medskip
\noindent
Now $F=E^\pi\cup D$. As all edges of $D$ has vertices with degree $1$ whereas all vertices of $e$ has degree $\ge 2$, we should have 
$$
e \in E^\pi
$$ which contradicts \eqref{condition 1}.

\medskip
\noindent
Hence we have $$\widetilde{\mathcal G}\not\cong
\widetilde{\mathcal H}$$
 
\medskip
\noindent
Now it remains to show that $$A_{\widetilde{\mathcal G}} = A_{\widetilde{\mathcal H}}.$$ 
Consider the following vertex ordering \[
\underbrace{u_1,\ldots,u_n}_{V},
\underbrace{a_{1,1},\ldots,a_{1,m-1},\ldots,
a_{m,1},\ldots,a_{m,m-1}}_{\text{additional vertices with degree }1}.
\] 
With respect to this ordering, the adjacency matrix of
$\widetilde{\mathcal G}$ has the block decomposition
\[
A_{\widetilde{\mathcal G}}
=
\begin{pmatrix}
A_{\mathcal G}+B& C\\
C^{T} & K
\end{pmatrix},
\]
where $A_{\mathcal G}$ is the adjacency matrix of the original hypergraph $\mathcal G$,  $C$ is the contribution of the adjacency relation between $S$ and $V$, $K$ records the adjacency relations within $S$, $B$ is the contribution of the added pendant hyperedges to
adjacencies among $V$.

\medskip
\noindent
Now observe that each pendant edges have only one unique vertices from $V$. So effectively $B=0$. Hence \[
A_{\widetilde{\mathcal G}}
=
\begin{pmatrix}
A_{\mathcal G}& C\\
C^{T} & K
\end{pmatrix}.
\]
Now observe that the added pendant hyperedges are identical in
$\widetilde{\mathcal G}$ and $\widetilde{\mathcal H}$,
the matrices $C$ and $K$ are same for both. Also we permuted the vertices from $V$ only. So
\[
A_{\widetilde{\mathcal G}^{\,\pi}}
=
\begin{pmatrix}
A_{\mathcal G^\pi} & C\\
C^{T} & K
\end{pmatrix}.
\]
As
$\pi\in\APP(\mathcal G)$,
definition~\eqref{def:adj-pres} gives
\[
A_{\mathcal G}
=
A_{\mathcal G^\pi}.
\]
Therefore
\[
A_{\widetilde{\mathcal G}}
=
A_{\widetilde{\mathcal H}}.
\]
\end{proof}

\noindent
We have shown that if we have a uniform hypergraph $\mathcal{G}$ with degree of all vertices at least $2$ and for which $\mathrm{APP}(\mathcal{G}) \setminus \mathrm{Aut}(\mathcal{G}) \neq \emptyset$,  we can construct two non-isomorphic coadjacent hypergraphs. Clealry coadjacent implies adjacency-cospectral. The natural question that arises now is whether there exists an infinite family of hypergraphs satisfying the above conditions. Otherwise, the proposed construction would become specific, limiting its broader applicability and significance. In the next section we consider hypergraphs coming form a certain class of projective planes and established that the above conditions hold for these hypergraphs.

\section{Hypergraphs From Finite Projective Planes}
\noindent
In this section we use some known structures and results from finite projective geometry to show existence of an infinite family of hypergraphs which satisfy all the conditions of Theorem \eqref{thm:final}. Before going into hypergraph let us define what a finite projective plane is.

\begin{definition}\cite{dougherty2020combinatorics}[Projective Plane]
A projective plane is defined as a triple $(\Pi,\Lambda,\sim)$, where
$\Pi$ denotes a set of points, $\Lambda$ denotes a set of lines, and
$\sim$ is an incidence relation between points and lines, satisfying
the following axioms:
\begin{enumerate}
    \item For every pair of distinct points $A,B\in\Pi$,
    there exists exactly one line $l\in\Lambda$ incident with both
    $A$ and $B$.

    \item For every pair of distinct lines $l_1,l_2\in\Lambda$,
    there exists exactly one point $C\in\Pi$ incident with both
    $l_1$ and $l_2$.

    \item There exist four points in $\Pi$, no three of which lie
    on the same line.
\end{enumerate}
\end{definition}

\noindent
This is a standard definition of a projective plane. We can think of it as a set of points, with lines defined as subsets of the set of points satisfying certain restrictions. In particular, any two distinct points lie on exactly one common line, and any two distinct lines have exactly one point in common. If both the set of points and the set of lines are finite, we call it a finite projective plane. We say that a finite projective plane has order $q$ if every line contains exactly $q+1$ points and every point lies on exactly $q+1$ lines. It is also a standard exercise to show that the numbers of points and lines in a finite projective plane of order $q$ are both equal to $q^2+q+1$; see \cite{dougherty2020combinatorics}. Here, we consider finite projective planes of order $q$, where $q\geq 2$. A very standard example is the Fano plane.
\begin{example}[The Fano plane]\label{ex:fano}
It is projective plane of order $2$. It has of $7$ points and $7$ lines. Each
line contains exactly $3$ points, every point lies on exactly $3$ lines, and
every pair of distinct points lies on a unique line.

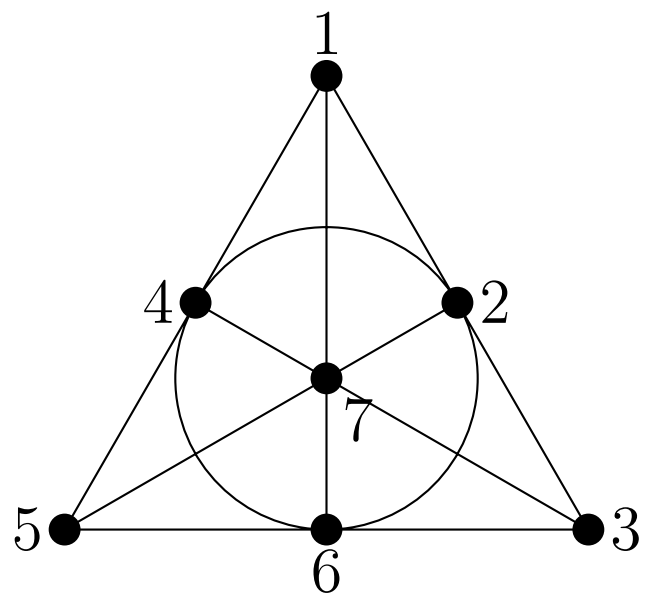
\begin{figure}[ht]
\centering
\begin{tikzpicture}
  \draw (30:1)  -- (210:2)
        (150:1) -- (330:2)
        (270:1) -- (90:2)
        (90:2)  -- (210:2) -- (330:2) -- cycle
        (0:0)   circle (1);

  \fill (0:0)   circle(3pt)
        (30:1)  circle(3pt)
        (90:2)  circle(3pt)
        (150:1) circle(3pt)
        (210:2) circle(3pt)
        (270:1) circle(3pt)
        (330:2) circle(3pt);

  \node[below]       at (0:0)   {$\quad7$};
  \node[right]       at (30:1)  {$2$};
  \node[above]       at (90:2)  {$1$};
  \node[left]        at (150:1) {$4$};
  \node[left]        at (210:2) {$5$};
  \node[below]       at (270:1) {$6$};
  \node[right]       at (330:2) {$3$};
\end{tikzpicture}
\caption{The Fano plane (finite projective plane of order \(2\)).}
\end{figure}

Its seven lines are
\[
\begin{aligned}
L_1&=\{1,2,3\}, &
L_2&=\{1,4,5\}, &
L_3&=\{1,6,7\},\\
L_4&=\{2,4,6\}, &
L_5&=\{2,5,7\}, &
L_6&=\{3,4,7\},\\
L_7&=\{3,5,6\}.
\end{aligned}
\]
\end{example}

\medskip
\noindent
A finite projective plane naturally gives rise to a hypergraph by taking its points as vertices and its lines as hyperedges. Connections between finite projective planes and hypergraph theory have been studied in the literature; see, for example, \cite{frankl1986finite,henning2018transversals}. Viewing the points as vertices and the lines as hyperedges gives the incidence hypergraph of the projective plane.

\begin{definition}[Incidence hypergraph of a finite projective plane]
Let $(P,L)$ be a finite projective plane of order $q$, where $P$ is the list of points and $L$ is the list of lines as a collection of $q+1$ subsets of $P$. The \emph{incidence hypergraph} of $(P,L)$ is the hypergraph $(V,E),$ whose vertex set $V=P$ and hyperedge set $E=L$. 

\medskip
\noindent
Observe from definition that incidence hypergraph coming from a projective plane of order $q$ is $q+1$ uniform and $q+1$ regular and linear. This is evident from the definition as each line contains $q+1$ points, each point incident to exactly $q+1$ lines and two line intersect exaclty at one point.
\end{definition}

\begin{example}\label{ex:4.4}
The projective plane in example \eqref{ex:fano}
gives rise to a \(3\)-uniform hypergraph $\mathcal{F}=(V,E),$ where $V=\{1,2,3,4,5,6,7\},$ and
\[
\begin{aligned}
E=\{&
\{1,2,3\},
\{1,4,5\},
\{1,6,7\},
\{2,4,6\},
\{2,5,7\},
\{3,4,7\},
\{3,5,6\}
\}.
\end{aligned}
\]

\noindent
Observe that if we delete the vertex labeled as $7$ and all incident edges then this is exactly the hypergraph in example \eqref{ex:coadjacent} up to a relabeling.
\end{example}

\medskip
\begin{proposition}\label{lem:adjmatrix-projective}
Let $\mathcal{G}$ be the hypergraph associated with a finite projective plane of order $q$. Let number of vertices of $\mathcal{G}$ be $n=q^2+q+1$. Then $A_{\mathcal G} = \frac1q(J_n-I_n),$ where \(I_n\) is the identity matrix of order \(n\), and \(J_n\) be the all-one matrix of order \(n\), i.e.,
\[
J_n = \mathbf{1}_n \mathbf{1}_n^T, \quad \mathbf{1}_n =
\begin{bmatrix}
1 & 1 & \cdots & 1
\end{bmatrix}^{T}
\in \mathbb{R}^n.
\]
\end{proposition}

\begin{proof}
If $i=j$, then by definition $(A_{\mathcal G})_{ii}=0.$ Suppose $i\neq j$. Since $\mathcal{G}$ is $(q+1)$-uniform and every pair of distinct vertices lies in exactly one common hyperedge, so by definition
\[
(A_{\mathcal G})_{ij}
=
\frac1{(q+1)-1}
=
\frac1q.
\]
Therefore, $$A_{\mathcal G}
=
\frac1q(J_n-I_n), \quad n=q^2+q+1.$$
\end{proof}

\medskip
\noindent
This is the main motivation behind chasing this type of hypergraphs. The reason we took this example because of its adjacency matrix, which is well behaved with any permutation of vertices. Consider the following corollary.
\begin{corollary}\label{thm:projective-app}
Let $\mathcal G$ be the hypergraph associated with a finite projective plane of order $q$. Then
\[
\APP(\mathcal G)=S_n, \quad n=q^2+q+1.
\]
  
\end{corollary}

\begin{proof}
By proposition~\eqref{lem:adjmatrix-projective}, 
$$A_{\mathcal G}
=
\frac1q(J_n-I_n).$$
Let $P$ be the permutation matrix corresponding to an arbitrary permutation $\pi\in S_n$. Since permutation matrices just permute rows and columns,
\[
PJP^T=J,
\qquad
PIP^T=I.
\]
Hence
\[
PA_{\mathcal G}P^T
=
\frac1q(PJP^T-PIP^T)
=
\frac1q(J-I)
=
A_{\mathcal G}.
\]
Since $A_{\mathcal G^\pi}
=
PA_{\mathcal G}P^T$ for every permutation $\pi\in S_n$, therefore every permutation is adjacency-preserving, and consequently
\[
\APP(\mathcal G)=S_n.
\]

\medskip
\noindent
There is also an interesting way to see this. Observe that the adjacency matrix for the hypergraph $\mathcal{G}$ is same as the adjacency matrix of the complete graph $K_n$ if we ignore the normalizing factor $\frac{1}{q}$. Also by Propesition \eqref{aut=app}, we know for simple graphs adjacency-preserving permutation group and automorphism groups are same. Hence, $$\pi \in \APP(\mathcal{G}) \iff \pi \in \Aut(K_n) =S_n.$$
\end{proof}

\medskip
\noindent
Clearly for the hypergraphs coming from finite projective planes of order $q\ge 2$, degree of every vertex is $q+1\ge 3$ as well as $q+1$ uniform. So it satisfies the first two condition of Theorem \eqref{thm:final}. Also we showed that for these type of hypergraphs the adjacency preserving group is exactly $S_n$. It remains to show that the automorphism group is strictly smaller that $S_n$.

\begin{proposition}\label{thm:strict-projective}
Let $\mathcal G = (V,E)$ be the incidence hypergraph associated with a finite projective plane of order $q\ge2$. Then
\[
\Aut(\mathcal G)\subsetneq\APP(\mathcal G).
\] 
\end{proposition}

\begin{proof}
By corollary~\eqref{thm:projective-app}, $\APP(\mathcal G)=S_n.$ Then it is enough to prove that $\Aut(\mathcal G)\neq S_n.$ Suppose for contradiction $$\Aut(\mathcal G)=S_n.$$

\noindent
Let $s\subset V$ with $|s|=q+1$. Take any $e \in E$. Clearly there exist a permuatation of the vertex set which maps $e$ to $s$. Now as any permutation of the vertices is an automorphism, $s$ is also and edge. Hence we get for any subset of $V$ with cardinality $q+1$ is an edge. So total number of distinct edges is $$\binom{q^2+q+1}{q+1}.$$

\noindent
But we know that number of edge for the hypergraph is $q^2+q+1$. Now 
$$\binom{q^2+q+1}{q+1} > q^2+q+1 , \quad q \ge 2 .$$

\noindent
Hence we get $\Aut(\mathcal G)\subsetneq S_n.$
\end{proof}
\begin{example}
Consider the hypergraph in example \eqref{ex:4.4}. By proposition~\eqref{lem:adjmatrix-projective},
\[
A_{\mathcal F}
=
\frac12(J_7-I_7).
\]
Therefore,
\[
\APP(\mathcal F)=S_7,
\]
while
\[
\Aut(\mathcal F)\subsetneq \APP(\mathcal F).
\]
Indeed, consider the transposition
\[
\pi=(1\;2).
\]
Since
\[
\APP(\mathcal F)=S_7,
\]
we have \(\pi\in\APP(\mathcal F)\). However,
\[
\pi(\{1,4,5\})=\{2,4,5\},
\]
which is not a hyperedge of \(\mathcal F\). Hence,
\[
\pi\notin\Aut(\mathcal F),
\]
\end{example}

\medskip
\noindent
Finally we have shown a class of hypergraph which satify all properties required to apply Theorem \eqref{thm:final}, that is hypergraphs coming from projective planes or order $q \ge 2$. It is a well-known result in finite geometry that, for every prime power $q=p^n$, there exists a finite projective plane of order $q$ \cite{schmidt2002characters}. Since there are infinitely many primes, it follows that there exist infinitely many such hypergraphs. 

\section{Conclusion}
\noindent
Our main idea behind this construction is basically generalizing the example \ref{ex:coadjacent-nonisomorphic}. Adjacency matrix may not always fully determine a hypergraph, we exploit this idea using the group difference $\APP(\mathcal G) \setminus \Aut(\mathcal G)$. Also the hypergraphs coming from finite projective planes serves as very good candidates for our construction. Observe that the non-isomorphic cospectral hypergraphs obtained in Example~\eqref{finalex} differ slightly from those presented in Example~\eqref{ex:coadjacent-nonisomorphic}. Indeed, the manner in which the new vertices and edges are attached is different, one is preserving linearity but use lots of extra vertices while other one loose linearity but use only one extra vertices.

\section*{Declarations}

\textbf{Competing interest :} There is no conflict of interest to declare.

\bibliographystyle{plain}
\bibliography{ref}

\end{document}